\documentclass[12pt]{amsart}
\usepackage{amsmath,amsthm,amssymb}
\usepackage[mathcal]{euscript}

\usepackage[all]{xy}

\newcommand{\smallxymatrix}[1]{\xymatrix@1@=1pc
                                       {#1}
                                   } 
\newcommand{\tinyxymatrix}[1]{
\def\objectstyle{\scriptstyle}
\def\labelstyle{\scriptstyle}
\vcenter{\xymatrix@-1.2pc{#1} }} 

\newcommand{\inpair}[2]{\ar@<.5ex>[r]^-{#1}\ar@<-.5ex>[r]_-{#2} } 
\newcommand{\smallpair}[4]{\smallxymatrix{ 
                      #1\ar@<.5ex>[r]^-{#2}\ar@<-.5ex>[r]_-{#3} & #4  
                                         }} 
\newcommand{\pair}[4]{\xymatrix@1{ 
                      #1\ar@<.5ex>[r]^-{#2}\ar@<-.5ex>[r]_-{#3} & #4  
                                      }} 
\newcommand{\xylabel}[2]{\ar@{}[#1]|{#2}} 
\newcommand{\xylabelc}[2]{\ar@{}[#1]|-{#2}} 
\newcommand{\xyinc}{\ar@{^{(}->}} 

\newcommand{\ten}{\mbox{\hspace*{-.5pt}\raisebox{1pt}{${\scriptstyle \otimes}$}
\hspace*{-4pt}}}

\newcommand{\st}{\mathrm{st}}

\newcommand{\odash}[2]{\mbox{$#1$\nobreakdash-\hspace{0pt}#2}} 

\newcommand{\infbi}{\odash{\epsilon}{bialgebra}}
\newcommand{\infbis}{\odash{\epsilon}{bialgebras}}
\newcommand{\quadri}{\text{quadri-algebra}}
\newcommand{\quadris}{\text{quadri-algebras}}
\newcommand{\quadrisub}{\text{quadri-subalgebra}}

\newcommand{\id}{\mathit{id}}

\newcommand{\End}{\textsf{End}}

\newcommand{\Sh}{\textsf{Sh}}     

\newcommand{\equal}[1]{
{\stackrel{{\textstyle #1}}{\ {\textstyle =}\ } }}
 
\newcommand{\map}[1]{\xrightarrow{#1}}

\renewcommand{\ne}{\nearrow} 
\newcommand{\se}{\searrow}
\newcommand{\sw}{\swarrow}
\newcommand{\nw}{\nwarrow}
\newcommand{\net}{\ne^t}
\newcommand{\set}{\se^t}
\newcommand{\swt}{\sw^t}
\newcommand{\nwt}{\nw^t}
\newcommand{\nep}{\ne^{op}}
\newcommand{\sep}{\se^{op}}
\newcommand{\swp}{\sw^{op}}
\newcommand{\nwp}{\nw^{op}}
\newcommand{\east}{\succ}
\newcommand{\west}{\prec}
\newcommand{\north}{\wedge}
\newcommand{\south}{\vee}
\newcommand{\eastt}{\east^t}
\newcommand{\westt}{\west^t}
\newcommand{\northt}{\north^t}
\newcommand{\southt}{\south^t}

\newcommand{\calQ}{\mathcal{Q}}    

\newcommand{\field}{\Bbbk} 

\theoremstyle{plain}
\newtheorem{theo}{Theorem}[section]
\newtheorem{prop}[theo]{Proposition}
\newtheorem{lemm}[theo]{Lemma}
\newtheorem{coro}[theo]{Corollary}

\theoremstyle{definition}
\newtheorem{defi}[theo]{Definition}
\newtheorem{exas}[theo]{Examples} 

\theoremstyle{remark}
\newtheorem{rems}[theo]{Remarks}  

\newcommand{\bd}{\begin{defi}}
\newcommand{\ed}{\end{defi}}
\newcommand{\bt}{\begin{theo}}
\newcommand{\et}{\end{theo}}
\newcommand{\bp}{\begin{prop}}
\newcommand{\ep}{\end{prop}}
\newcommand{\bl}{\begin{lemm}}
\newcommand{\el}{\end{lemm}}
\newcommand{\bc}{\begin{coro}}
\newcommand{\ec}{\end{coro}}
\newcommand{\bpf}{\begin{proof}}
\newcommand{\epf}{\end{proof}}
\newcommand{\br}{\begin{rems}\begin{flushleft}\end{flushleft}\nopagebreak} 
\newcommand{\er}{\end{rems}}
\newcommand{\bex}{\begin{exas}\begin{flushleft}\end{flushleft}\nopagebreak} 
\newcommand{\eex}{\end{exas}}

\newcommand{\be}{\begin{enumerate}}
\newcommand{\ee}{\end{enumerate}}
\newcommand{\bi}{\begin{itemize}}
\newcommand{\ei}{\end{itemize}}

\begin{document}

\title{Quadri-algebras}
\author[M. Aguiar]{Marcelo Aguiar}
\address{Department of Mathematics\\
         Texas A\&M University\\
         College Station\\
         TX  77843\\
         USA}
\email{maguiar@math.tamu.edu}
\urladdr{http://www.math.tamu.edu/$\sim${}maguiar}

\author[J.-L. Loday]{Jean-Louis Loday}
\address{Institut de Recherche Math\'ematique Avanc\'ee\\
CNRS et Universit\'e Louis Pasteur\\
7 rue R. Descartes\\
67084 Strasbourg Cedex, France}
\email{loday@math.u-strasbg.fr}
\urladdr{http://www-irma.u-strasbg.fr/$\sim${}loday/}

\thanks{The second author would like to thank Texas A\&M University for an invitation in October
2002 during which this work was started.}

\date{September 9, 2003}

\begin{abstract} We introduce the notion of $\quadris$. 
These are associative algebras for which the multiplication can be decomposed as the sum
of four operations  in a certain coherent manner. We present several examples of $\quadris$:
the algebra of permutations, the shuffle algebra, tensor products of dendriform algebras.
We show that a pair of commuting
Baxter operators on an associative algebra gives rise to a canonical $\quadri$ structure
on the underlying space of the algebra. The main example is provided by the algebra
$\End(A)$ of linear endomorphisms of an infinitesimal bialgebra $A$. This algebra
carries a canonical pair of commuting Baxter operators: $\beta(T)=T\ast\id$ and
$\gamma(T)=\id\ast T$, where
$\ast$ denotes the convolution of endomorphisms. It follows that $\End(A)$ is a
$\quadri$, whenever $A$ is an infinitesimal bialgebra. We also discuss commutative
$\quadris$ and state some conjectures on the free $\quadri$.
\end{abstract}
\maketitle

\section*{Introduction} \label{S:int}

The study of the space of endomorphisms of an infinitesimal bialgebra 
 revealed the existence of peculiar algebraic structures. More 
 specifically, the convolution of endomorphisms gives rise to a pair of 
 commuting Baxter operators 
\[ \beta(T)=T\ast\id \text{ \ and \ }\gamma(T)=\id\ast T\]
and each of these determines a 
dendriform  structure on the space of endomorphisms. Obviously these two 
 structures  are somehow intertwined, but under which rule? 

In this work we answer this question. We introduce the notion of 
 ``quadri-algebras". A quadri-algebra 
 is an associative algebra whose multiplication is the sum of 4 
 operations. These operations satisfy 9 relations. We show that the 
 space of endomorphisms of an infinitesimal bialgebra has a natural 
 structure of quadri-algebra which encompasses the two dendriform structures. 

 In Section~\ref{S:quadri} we recall the definition of dendriform algebras and we 
 introduce quadri-algebras. We also provide a number of examples of quadri-algebras and
constructions relating quadri-algebras to dendriform algebras. We also discuss
``commutative'' quadri-algebras.

  In Section~\ref{S:baxter} we recall the definition of Baxter operators and we 
 show that a pair of commuting Baxter operators gives rise to a 
 quadri-algebra. 

 In Section~\ref{S:infbi} we apply the preceding result to show that the 
space of endomorphisms of an infinitesimal bialgebra may be naturally 
 endowed with a  quadri-algebra structure. 

  In Section~\ref{S:free} we conclude with some conjectures related to the free
quadri-algebra  on one generator. In particular, we conjecture that the dimension \
of its homogeneous component of degree $n$ is equal to the number of   non-crossing
 connected graphs on $n+1$ vertices.

\section{Quadri-algebras}\label{S:quadri}

\subsection{Dendriform algebras} A {\em dendriform algebra} is a vector space $D$
together with two operations
 $\west :D\ten D\to D$ and $\east :D\ten D\to D$, called  {\it left} and {\it right} respectively, such
that
\begin{align}
(x\west y)\west z &= x\west (y\west z)+x\west (y\east z) \notag\\
(x\east y)\west z &= x\east (y\west z) \label{E:dendri}\\
(x\west y)\east z+(x\east y)\east z & = x\east (y\east z)\,. \notag
\end{align}

Dendriform algebras were introduced by the second author~\cite[Chapter 5]{L}.
See~\cite{R,LR1,LR2,Ler,depaul} for additional work on this subject.

Defining a new operation by
\begin{equation}\label{E:dendri-assoc}
x\star y:=x\west y+x\east y
\end{equation}
 permits us to rewrite axioms~\eqref{E:dendri} as 
\begin{align*}
(x\west y)\west z &= x\west (y\star z) \notag\\
(x\east y)\west z &= x\east (y\west z) \\
(x\star y)\east z & = x\east (y\east z)\,. \notag
\end{align*}

By adding the three relations we see that the operation $\star$ is associative.  For this reason,
a dendriform algebra may be regarded as  an associative algebra $(D,\star)$ for which the
multiplication $\star$ can be decomposed as the sum of two coherent operations.

We now introduce a class of algebras with
an associative multiplication which can be coherently decomposed as the sum of four
operations. As a consequence, these algebras carry two distinct dendriform  structures. 

\subsection{Quadri-algebras. Definition}\label{S:quadri-def} A {\em quadri-algebra} is a vector space $Q$
together with four operations $\se$, $\ne$,
$\nw$ and $\sw:Q\ten Q\to Q$ satisfying the axioms below~\eqref{E:quadri}. In
order to state them, consider the following operations:
\begin{align}
x\east y:= &x\ne y+x\se y \label{E:east}\\
x\west y:= &x\nw y+x\sw y \label{E:west}\\
x\south y:= &x\se y+x\sw y \label{E:north}\\
x\north y:= &x\ne y+x\nw y \label{E:south}\\
\intertext{and}
x\star y:=&x\se y+x\ne y+x\nw y+x\sw y \label{E:star}\\= & x\east y+x\west y=x\south y+x\north y\,. \notag
\end{align}

The axioms are
\begin{equation}  \label{E:quadri}
\raisebox{25pt}{\xymatrix@C=5pt@R=5pt{
{\mbox{\footnotesize $(x\nw y)\nw z =x\nw(y\star z)$}} & &
{\mbox{\footnotesize $(x\ne y)\nw z = x\ne(y\west z)$}} &  &
{\mbox{\footnotesize $(x\north y)\ne z=x\ne(y\east z)$}}  \\
 {\mbox{\footnotesize $(x\sw y)\nw z =x\sw(y\north z)$}} & &
 {\mbox{\footnotesize $(x\se y)\nw z = x\se(y\nw z)$}} &  &
{\mbox{\footnotesize $(x\south y)\ne z=x\se(y\ne z)$}}             \\
 {\mbox{\footnotesize $(x\west y)\sw z=x\sw(y\south z)$}} & &
 {\mbox{\footnotesize $(x\east y)\sw z = x\se(y\sw z)$}} &  &
{\mbox{\footnotesize $(x\star y)\se z =x\se(y\se z)$}}
}}
\end{equation}
\smallskip

 We refer to the operations $\se$, $\ne$, $\nw$, $\sw$,
as {\em southeast, northeast, northwest,} and {\em southwest}, respectively. Accordingly,
we use {\em north, south, west,} and {\em east} for $\north$, $\south$, $\west$ and $\east$.

 The axioms are displayed in the form of a
$3\times 3$ matrix. We will make use of standard matrix terminology (entries, rows and
columns) to refer to them.

\subsection{From quadri-algebras to dendriform algebras}\label{S:quadri-dendri}

The three column sums of the matrix of axioms~\eqref{E:quadri} yield
\[(x\west y)\west z=x\west(y\star z), \ \ (x\east y)\west z = x\east(y\west z) \text{
\ and \ } (x\star y)\east z =x\east(y\east z)\,.\]
Thus, endowed with the operations west for left and east for right,
$Q$ is a dendriform algebra. We denote it by $Q_h$ and call it the {\em horizontal} dendriform
algebra associated to $Q$.

Consider instead the three row sums in~\eqref{E:quadri}. We obtain
\[(x\north y)\north z=x\north(y\star z), \ \  (x\south y)\north z = x\south(y\north z)
\text{ \ and \ }
(x\star y)\south z =x\south(y\south z)\,.\]
Thus, endowed with the operations north for left  and south for right,
$Q$ is a dendriform algebra.  We denote it by $Q_v$ and call it the {\em vertical} dendriform
algebra associated to $Q$.

The associative operations corresponding to the dendriform algebras $Q_h$ and $Q_v$ by means
of~\eqref{E:dendri-assoc} coincide, according to~\eqref{E:star}. Thus, a $\quadri$ may
be seen as an associative algebra $(Q,\star)$ for which the multiplication $\star$ can be decomposed into four coherent operations.

\subsection{Example. The algebra of permutations}\label{S:permutations} Let $S_n$ denote the symmetric group on $n$
letters and denote by $\field S_n$ the vector space spanned by its elements (over a field
$\field$). Consider the spaces
\[\field S_\infty:=\bigoplus_{n\geq 0}\field S_n \text{ \ and \ }\field S_{\geq
2}:=\bigoplus_{n\geq 2}\field S_n\,.\]
 Let $\Sh(p,q)$ denote the set of \odash{(p,q)}{shuffles}, that is, those
permutations $\zeta\in S_{p+q}$ such that
\[\zeta(1)<\ldots<\zeta(p) \text{ \ and \ }  \zeta(p+1)<\ldots<\zeta(p+q)\,.\]
Note that any such shuffle satisfies
\[\zeta^{-1}(1)=1 \text{ or }p+1 \text{ \ \ and \ \ } \zeta^{-1}(p+q)=p \text{ or }p+q\,.\]
Therefore, if both $p$ and $q$ are at least $2$, the set of  \odash{(p,q)}{shuffles}
decomposes into the following four disjoint subsets:
\begin{align*}
\Sh^1(p,q) &=\{\zeta\in \Sh(p,q)\ \mid\  \zeta^{-1}(1)=1  \text{ and } \zeta^{-1}(p+q)=p+q\}\\
\Sh^2(p,q) &=\{\zeta\in \Sh(p,q)\ \mid\  \zeta^{-1}(1)=p+1  \text{ and } \zeta^{-1}(p+q)=p+q\}\\
\Sh^3(p,q) &=\{\zeta\in \Sh(p,q)\ \mid\  \zeta^{-1}(1)=1  \text{ and } \zeta^{-1}(p+q)=p\}\\
\Sh^4(p,q) &=\{\zeta\in \Sh(p,q)\ \mid\  \zeta^{-1}(1)=p+1  \text{ and } \zeta^{-1}(p+q)=p\}\,.
\end{align*}
For $\sigma\in S_p$ and $\tau\in S_q$, let $\sigma\times\tau$ denote the
permutation in $S_{p+q}$ defined by
\[i\mapsto\begin{cases} \sigma(i) & \text{ if }i\leq p\,,\\
p+\tau(i-p) & \text{ if }i>p\,. \end{cases}\]
The composition of permutations $\zeta,\xi\in S_n$ is $(\zeta\cdot\xi)(i)=\zeta\bigl(\xi(i)\bigr)$.
Let $\sigma$ and $\tau$ be as above, with $p,q\geq 2$. Define operations on $\field S_{\geq
2}$ by
\begin{align*}
\sigma \sw \tau &= \sum_{\zeta\in\Sh^{1}(p,q)}\zeta\cdot(\sigma\times\tau) \\
\sigma \se \tau &= \sum_{\zeta\in\Sh^{2}(p,q)}\zeta\cdot(\sigma\times\tau) \\
\sigma \nw \tau &=  \sum_{\zeta\in\Sh^{3}(p,q)}\zeta\cdot(\sigma\times\tau)\\
\sigma \ne \tau &= \sum_{\zeta\in\Sh^{4}(p,q)}\zeta\cdot(\sigma\times\tau)\,. \\
\end{align*}
This endows $\field S_{\geq 2}$ with a $\quadri$ structure. Axioms~\eqref{E:quadri} are
easily verified with the aid of the standard bijection
\[\Sh(p+q,r)\times\Sh(p,q)\cong \Sh(p,q+r)\times\Sh(q,r)\,.\]
The vertical dendriform structure corresponding to this $\quadri$ is the dendriform
algebra of permutations introduced in~\cite[Section 5.4]{L} and 
further studied in~\cite{LR2} (this structure is in fact defined on $\field S_{\geq 1}$).

The associative structure is
\[\sigma \star \tau = \sum_{\zeta\in\Sh(p,q)}\zeta\cdot(\sigma\times\tau)\,.\]
This is the multiplication of Malvenuto and Reutenauer~\cite[pages 977-978]{MR} (which is
in fact defined on the whole $\field S_\infty$).

\subsection{Tensor product of dendriform algebras}\label{S:tensor}
The tensor product of two dendriform algebras carries a natural
$\quadri$ structure.

Let $A$ and $B$ be two dendriform algebras. On the space $A\ten B$, define operations
\begin{align*}
(a_1\ten b_1)\nw(a_2\ten b_2) &=(a_1\west a_2)\ten(b_1\west b_2)\\
(a_1\ten b_1)\sw(a_2\ten b_2) &=(a_1\west a_2)\ten(b_1\east b_2)\\
(a_1\ten b_1)\ne(a_2\ten b_2) &=(a_1\east a_2)\ten(b_1\west b_2)\\
(a_1\ten b_1)\se(a_2\ten b_2) &=(a_1\east a_2)\ten(b_1\east b_2)\,.
\end{align*}
Axioms~\eqref{E:quadri} follow from axioms~\eqref{E:dendri} for $A$ and $B$.

The corresponding horizontal dendriform structure is
\begin{align*}
(a_1\ten b_1)\west(a_2\ten b_2) &=(a_1\west a_2)\ten(b_1\star b_2)\\
(a_1\ten b_1)\east(a_2\ten b_2) &=(a_1\east a_2)\ten(b_1\star b_2)\,.
\end{align*}
This is the structure used in the definition of dendriform bialgebras~\cite{R}.

The corresponding associative structure is simply the tensor product of the
associative structures on $A$ and $B$:
\[(a_1\ten b_1)\star(a_2\ten b_2) =(a_1\star a_2)\ten(b_1\star b_2)\,.\]

\subsection{Opposite and transpose of a quadri-algebra}\label{S:opptrans}

The matrix of axioms~\eqref{E:quadri} is symmetric with respect to the
main diagonal, provided we interchange $\sw$ with $\ne$, and leave the other two
operations unchanged (note this has the effect of interchanging $\north$ with $\west$
and $\south$ with $\east$). Therefore, starting from a $\quadri$ $(Q,\se,\ne,\nw,\sw)$
and defining
\[x\set  y=x\se y,\ \ \ \ x\net  y=x\sw y,\ \ \ \ x\nwt  y=x\nw y
\text{ \ \ \ and \ \ \ }x\swt  y=x\ne y\,,\]
 one obtains a  new $\quadri$ structure on the underlying space of $Q$. We refer to this
new $\quadri$ as the {\em transpose} of $Q$ and denote it by $Q^t$. The dendriform structures corresponding to $Q$ (as in~\ref{S:quadri-dendri})
are the same as those for $Q^t$, in the opposite order: $Q_h=Q^t_v$ and $Q_v=Q^t_h$.
In fact,
\begin{gather*}
x\southt y=x\set  y+x\swt y =x\se y+x\ne y=x\east y\\
x\northt y=x\net  y+x\nwt y =x\sw y+x\nw y=x\west y\\
x\eastt y=x\net  y+x\set y =x\sw y+x\se y=x\south y\\
x\westt y=x\nwt  y+x\swt y =x\nw y+x\ne y=x\north y\\
\end{gather*}
It follows that the associative structures corresponding to $Q$ and $Q^t$ coincide.

The matrix of axioms exhibits another symmetry, with respect to the center (the entry $(2,2)$),
provided we replace each arrow by its
opposite ($\nw$ by $\se$ and $\sw$ by $\ne$) and reverse the order of the variables.
Therefore, starting from a $\quadri$ $(Q,\se,\ne,\nw,\sw)$
and defining
\[x\sep  y=y\nw x,\ \ \  x\nep  y=y\sw x,\ \ \  x\nwp  y=y\se x
\text{ \ \ \ and \ \ \ }x\swp  y=y\ne x\,,\]
 one obtains a  new $\quadri$ structure on the underlying space of $Q$. We refer to this
new $\quadri$ as the {\em opposite} of $Q$ and denote it by $Q^{op}$.

\subsection{Commutative quadri-algebras}\label{S:commquadri}
A $\quadri$ $Q$ is said to be {\em commutative} if it coincides with its opposite: $Q=Q^{op}$.
Explicitly, this means that
\[ x\se  y=y\nw x \text{ \ and \ } x\ne  y=y\sw x\,.\]
One may then restate axioms~\eqref{E:quadri} in terms of the operations $\se$ and $\sw$ only.
One obtains that a commutative
$\quadri$ is a space $Q$ equipped with two operations $\se$ and $\sw:Q\ten Q\to Q$ such that
\begin{align}
x\se(y\se z) &=(x\star y)\se z \notag\\
x\se(y\sw z) &=(x\east y)\sw z \label{E:commquadri} \\
x\se(y\sw z) &=y\sw(x\south z) \notag\end{align}
where we have set, in agreement with the notation for general $\quadris$,
\begin{align*}
x\east y:= &x\se y+y\sw x \\
x\south y:= &x\se y+x\sw y \\
\intertext{and}
x\star y:=&x\se y+x\sw y+y\se x+y\sw x \\
= & x\east y+y\east x=x\south y+y\south x \,.
\end{align*}

It follows from~\eqref{E:commquadri} that the operation $\star$ is associative and
commutative, and that each of the operations $\east$ and $\south$ define a left
Zinbiel structure on $A$ (see~\cite[Section 7]{L} for the definition of Zinbiel algebras).

\subsection{Example. The shuffle algebra}\label{Ex:shuffle}
The shuffle algebra of a vector space $V$ provides an example of a commutative $\quadri$.

Consider the vector spaces
\[T(V) := \bigoplus_{n\geq 0} V^{\ten n} \text{ \ and \ }
T^{\geq 2}(V) := \bigoplus_{n\geq 2} V^{\ten n}\,.\]
We adopt the concatenation notation for elements in $V^{\ten n}$.
In the sequel $a,b,c,d \in V$ and $\omega, \theta\in V^{\ten n}$ for some $n\geq 0$.
On the space $ T^{\geq 2}(V)$ define operations as follows:

\begin{align*}
a\omega b \sw c\theta d &= a(\omega b \star c \theta )d \\
a\omega b \se c\theta d &= c(a\omega b \star \theta )d \\
a\omega b \nw c\theta d &= a(\omega \star c \theta d)b \\
a\omega b \ne c\theta d &= c(a\omega  \star  \theta d)b\,. \\
\end{align*}where, as before, $x\star y:=x\se y+x\ne y+x\nw y+x\sw y$. In low dimensions we start with
$a\star b = ab + ba$ and $1$ is a unit for $\star\ $ .  It follows that $\star$ is the
shuffle product (which is defined on the whole space $T(V)$) and
\begin{align*}
a\omega b \west c\theta d &= a(\omega b \star c \theta d) \\
a\omega b \east c\theta d &= c(a\omega b \star \theta d) \\
a\omega b \south c\theta d &= (a\omega  \star c \theta d)b \\
a\omega b \north c\theta d &= (a\omega b \star c \theta )d \,,
\end{align*}

Axioms~\eqref{E:quadri} can now be easily verified. For instance,
\[(a\omega b\nw c\theta d)\nw e\xi f = a(\omega \star c\theta d)b \nw e\xi f =
a(\omega \star c\theta d\star e\xi f)b\,,\]
while
\[a\omega b \nw (c\theta d \star e\xi f) = a\omega b \nw (c\theta d \star e\xi f) =
a(\omega \star c\theta d\star e\xi f)b\,.\]
Thus, the axiom in entry $(1,1)$ holds by associativity of the shuffle product.
The other axioms can be verified similarly. In fact the 9 monomials obtained by formulas (8) begin with
$a$ (resp. $c$, resp. $e$) in the first (resp. second, resp. third) column and end   with
$b$ (resp. $d$, resp. $f$) in the first (resp. second, resp. third) row. 

This is actually an example of a commutative $\quadri$. In fact,
\[a\omega b \se c\theta d=c(a\omega b \star \theta )d= c(\theta\star a\omega b)d=
c\theta d\nw a\omega b\,,\]
by the commutativity of the shuffle product. Similarly, $a\omega b \ne c\theta d=c\theta d\sw a\omega b$.

Observe that the dendriform structures are defined on $T^{\geq 1}(V)$ and the associative 
structure is defined on $T(V)$ (the shuffle algebra).

\subsection{A generalization of quadri-algebras}\label{S:ennea}

There is a generalization  of dendriform algebras called {\em dendriform trialgebras}; they
carry $3$ operations satisfying  $7$ relations~\cite{LR3}. The free  dendriform trialgebra
can be explicitly described by means of  planar  rooted trees (not necessarily binary). It is
natural to expect that the notion of quadri-algebra has a  similar generalization involving
$3^2= 9$ operations satisfying $7^2 = 49$  relations. This has been found recently by
Leroux~\cite{Ler2}.

\section{Baxter operators}\label{S:baxter}

\subsection{Baxter operators on associative algebras}\label{S:baxter-assoc}
 Let $(A,\cdot)$ be an associative algebra. A {\em Baxter operator} is a map
$\beta:A\to A$ such that
\begin{equation}\label{E:baxter-assoc}
\beta(x)\cdot\beta(y)=\beta\Bigl(x\cdot\beta(y)+\beta(x)\cdot y\Bigr)\,.
\end{equation}
This identity appeared originally in the work of Glen Baxter~\cite{Bax}.
The importance of Baxter operators was emphasized by Gian-Carlo Rota~\cite{R1,R2}.

Given such $A$ and $\beta$, one may define new operations on $A$ by
\[x\west_\beta y=x\cdot\beta(y) \text{ \ and \ } x\east_\beta y=\beta(x)\cdot y\,.\]
It is easy to see that then $(A,\west_\beta ,\east_\beta )$ is a dendriform
algebra~\cite[Proposition 4.5]{depaul}. The resulting associative structure
$\star_\beta$ on $A$ is related to the original one as follows:
$\beta:(A,\star_\beta)\to (A,\cdot)$ is a morphism of associative algebras.

Similarly, starting from a Baxter operator on a dendriform algebra, one can construct
a $\quadri$ structure on the same space, as we explain next.

\subsection{Baxter operators on dendriform algebras}\label{S:baxter-dendri}
Let $(D,\west,\east)$ be a dendriform algebra. A Baxter
operator is a map $\gamma:D\to D$ such that
\begin{align}
\gamma(x)\east\gamma(y)=\gamma\Bigl(x\east\gamma(y)+\gamma(x)\east
y\Bigr)\,,\label{E:baxter-east}\\
\gamma(x)\west\gamma(y)=\gamma\Bigl(x\west\gamma(y)+\gamma(x)\west
y\Bigr)\,.\label{E:baxter-west}
\end{align}

Adding these equations we see that $\gamma$ is also a Baxter operator on the
associative algebra $(D,\star)$:
\begin{equation}\label{E:gamma-star}
\gamma(x)\star\gamma(y)=\gamma\Bigl(x\star\gamma(y)+\gamma(x)\star
y\Bigr)\,.
\end{equation}

\subsection{Proposition}\label{P:baxter-dendri}
{\em Let $(D,\west,\east)$ be a dendriform algebra and $\gamma:D\to D$ a Baxter operator.
 Define  new operations on $A$ by
\[x\se_\gamma y=\gamma(x) \east y,\ \  x\ne_\gamma y=x\east\gamma(y),\ \  x\sw_\gamma
y=\gamma(x)\west y
 \text{ \ \ and \ \ } x\nw_\gamma y=x\west \gamma(y)\,.\]
Then $(D,\se_\gamma ,\ne_\gamma,\nw_\gamma,\sw_\gamma )$ is a $\quadri$.}
\bpf  We verify  the axioms corresponding to the entries $(1,1)$ and $(2,2)$
in~\eqref{E:quadri}; the others are similar. We have
\begin{align*}
(x\nw_\gamma y)\nw_\gamma z&=\bigl(x\west \gamma(y)\bigr)\west\gamma(z)\equal{\eqref{E:dendri}}
x\west\bigl(\gamma(y)\star\gamma(z)\bigr)\\
&\equal{\eqref{E:gamma-star}}x\west\gamma\Bigl(x\star\gamma(y)+\gamma(x)\star y\Bigr)=x\west\gamma(y\star_\gamma z)=
x\nw_\gamma(y\star_\gamma z)\,.
\end{align*}
 Also,
\begin{align*}
(x\se_\gamma y)\nw_\gamma z&=\bigl(\gamma(x)\east y\bigr)\west\gamma(z)\\
&\equal{\eqref{E:dendri}}\gamma(x)\east\bigl(y\west\gamma(z)\bigr)=
x\se_\gamma(y\nw_\gamma z)\,,
\end{align*}
as needed.
\epf

Let $Q$ denote the resulting $\quadri$. The horizontal dendriform structure associated to $Q$ is
\begin{align*}
x\west_\gamma y=x\nw_\gamma y+x\sw_\gamma y=x\west\gamma(y)+\gamma(x)\west y\,,\\
x\east_\gamma y=x\se_\gamma y+x\ne_\gamma y=\gamma(x)\east y+x\east\gamma(y)\,.
\end{align*}
Therefore, axioms~\eqref{E:baxter-east}
and~\eqref{E:baxter-west}  can be rewritten as follows:
\[ \gamma(x)\east\gamma(y)=\gamma(x\east_\gamma y) \text{ \ and \ }
\gamma(x)\west\gamma(y)=\gamma(x\west_\gamma y)\,.\]
Thus, $\gamma$ is a morphism of dendriform algebras $Q_h\to D$.

On the other hand, the vertical dendriform structure associated to $Q$ is
\begin{align*}
x\north_\gamma y=x\ne_\gamma y+x\nw_\gamma y=x\east\gamma(y)+x\west\gamma(y)=x\star\gamma(y)\,,\\
x\south_\gamma y=x\se_\gamma y+x\sw_\gamma y=\gamma(x)\east y+\gamma(x)\west y=\gamma(x)\star y\,.
\end{align*}
Thus, $Q_v$ is the dendriform structure corresponding to $\gamma$ viewed as a Baxter
operator on the associative algebra $(D,\star)$.

From either of the two previous remarks it follows that $\gamma$ is a morphism of
associative algebras $(Q,\star_\gamma)\to (D,\star)$ (a fact that was implicitly used
in the proof of~\ref{P:baxter-dendri}).

\subsection{Pairs of commuting Baxter operators}\label{S:baxterpair}
Let $A$ be an associative algebra and $\beta$ and $\gamma$ two Baxter operators on $A$ that
commute, that is,
\[\beta\gamma=\gamma\beta\,.\]
In this situation, it is possible to construct a $\quadri$ structure on the underlying
space of $A$, by successively applying the constructions of~\ref{S:baxter-assoc}
and~\ref{P:baxter-dendri}, as we explain next.

\subsection{Proposition}\label{P:baxterpair}
{\em Let $\beta$ and $\gamma$ be a pair of commuting Baxter operators
on an associative algebra $A$.   Then $\gamma$ is
a Baxter operator on the dendriform algebra $(A,\west_\beta,\east_\beta)$
corresponding to $\beta$ as in~\ref{S:baxter-assoc}.}
\bpf  We verify axiom~\eqref{E:baxter-east}; axiom~\eqref{E:baxter-west} is similar.
\begin{align*}
 \gamma(x)\east_\beta\gamma(y) &=\beta\bigl(\gamma(x)\bigr)\cdot \gamma(y)
 =\gamma\bigl(\beta(x)\bigr)\cdot \gamma(y)\\
 &\equal{\eqref{E:baxter-assoc}}\gamma\Bigl(\beta(x)\cdot\gamma(y)+ \gamma\bigl(\beta(x)\bigr)
 \cdot y\Bigr)=\gamma\Bigl(\beta(x)\cdot\gamma(y)+ \beta\bigl(\gamma(x)\bigr)
 \cdot y\Bigr)\\
 &=\gamma\bigl(x\east_\beta\gamma(y)+\gamma(x)\east_\beta y \bigr)\,.
\end{align*}\epf

\subsection{Corollary}\label{C:baxterpair}
{\em Let $\beta$ and $\gamma$ be a pair of commuting Baxter operators
on an associative algebra $A$. Then there is a $\quadri$ structure on the
underlying space of $A$, with operations defined by }
\begin{align*}
x\se y&=\beta\bigl(\gamma(x)\bigr)y=\gamma\bigl(\beta(x)\bigr)y\,,\\
x\ne y &=\beta(x)\gamma(y)\,,\\
x\sw y &=\gamma(x)\beta(y)\,,\\
x\nw y &=x\beta\bigl(\gamma(y)\bigr)=x\gamma\bigl(\beta(y)\bigr)\,.
\end{align*}
\bpf
Apply Proposition~\ref{P:baxter-dendri} to the Baxter operator $\gamma$ on
the dendriform algebra $(A,\west_\beta,\east_\beta)$.
\epf

One may start by first constructing the dendriform algebra $(A,\west_\gamma,\east_\gamma)$
instead. By Proposition~\ref{P:baxterpair}, $\beta$ is a Baxter operator on this dendriform
algebra. Hence, Proposition~\ref{P:baxter-dendri} yields a new $\quadri$ structure on the
underlying space of $A$. It is easy to see that this is the transpose of the structure of
Corollary~\ref{C:baxterpair}, in the sense of~\ref{S:opptrans}.

Let $(A,\star)$ denote the associative structure corresponding to the $\quadri$ structure
on $A$ (according to~\ref{S:opptrans}, this is the same for one structure and its transpose). It follows from the remarks in~\ref{S:baxter-dendri} that both $\gamma$ and $\beta$ are
morphisms of associative algebras $(A,\star)\to (A,\cdot)$. In this sense, the pair of
commuting Baxter operators $\beta$ and $\gamma$ breaks the associativity of $A$ into
four pieces.

\subsection{Example}\label{Ex:AYB} Let $A$ be an associative algebra
and let $r=\sum_i u_i\ten v_i\in A\ten A$ be a solution of the associative Yang-Baxter
equation~\cite{cor,analog}. The map $\beta_r:A\to A$ given by 
\[\beta_r(x):=\sum_i u_ixv_i\]
is a Baxter operator on $A$~\cite[Section 5]{poisson}. Suppose that 
$s\in A\ten A$ is another solution of the associative Yang-Baxter equation.
If $r$ and $s$ commute as elements of $A\ten A^{op}$ then the operators
$\beta_r$ and $\beta_s$ commute and Corollary~\ref{C:baxterpair} applies.

\section{Quadri-algebras from infinitesimal bialgebras}\label{S:infbi}

\subsection{Infinitesimal bialgebras} An {\em infinitesimal bialgebra} (abbreviated $\infbi$) is
a triple $(A, \mu,\Delta)$ where
$(A, \mu)$  is an  associative algebra,  $(A,\Delta)$  is a coassociative coalgebra,
and for each $a,b\in A$,
\begin{equation}
\Delta(ab)= ab_1\ten b_2+ a_1\ten a_2b\,.\label{E:infbi}
\end{equation}
We write $\Delta(a)=a_1\ten a_2$ (simplified Sweedler's notation) and
$ (\Delta\ten \id )\Delta(a)=a_1\ten a_2\ten a_3$.
Infinitesimal bialgebras originated in the work of Joni and Rota~\cite{JR}.
See~\cite{cor,analog,depaul} for the basic theory of $\infbis$.

We view the space $\End(A)$ of linear endomorphisms of $A$ as an associative algebra
under composition, denoted simply by concatenation:
\[TS:A\map{S}A\map{T}A\,.\]
We make use of a second associative product on $\End(A)$, the convolution of endomorphisms,
defined by
\[T\ast S:A\map{\Delta}A\ten A\map{T\ten S}A\ten A\map{\mu}A\,.\]
Note that the associativity of this product does not depend on any compatibility
condition between $\Delta$ and $\mu$.

\subsection{Proposition}\label{P:baxter-end}
{\em Let $A$ be an $\infbi$ and consider $\End(A)$ as an associative algebra under composition.
There is a pair of commuting Baxter operators $\beta$ and $\gamma$ on $\End(A)$, defined by }
\[\beta(T)=\id\ast T \text{ and }\gamma(T)=T\ast\id\,.\]
\bpf The operators commute because the convolution product is associative. Let us check
axiom~\eqref{E:baxter-assoc} for the operator $\beta$; the verification for $\gamma$ is similar.
We have $\beta(S)(a)=a_1S(a_2)$. Hence, by~\eqref{E:infbi},
\[\Delta\bigl(\beta(S)(a)\bigr)=a_1S(a_2)_1\ten S(a_2)_2+a_1\ten a_2S(a_3)\,.\]
Therefore,
\begin{align*}
\beta(T)\beta(S)(a)&=a_1S(a_2)_1T\bigl(S(a_2)_2\bigr)+a_1T\bigl(a_2S(a_3)\bigr) \\
&=a_1\bigl(\beta(T)S\bigr)(a_2)+a_1\bigl(T\beta(S)\bigr)(a_2)\\
&=\beta\bigl(\beta(T)S\bigr)(a)+\beta\bigl(T\beta(S)\bigr)(a)\,,
\end{align*}
as needed.
\epf

\subsection{Corollary}\label{C:quadri-end}
{\em Let $A$ be an $\infbi$. There is a $\quadri$
structure on the space $\End(A)$ defined by }
\begin{align*}
T\se S&=(\id\ast T\ast \id)S\,,\\
T\ne S &=(\id\ast T)(S\ast\id)\,,\\
T\sw S &=(T\ast\id)(\id\ast S)\,,\\
T\nw S &=T(\id\ast S\ast \id)\,.
\end{align*}
\bpf
Apply Corollary~\ref{C:baxterpair} to the pair of commuting Baxter operators of
Proposition~\ref{P:baxter-end}.
\epf

The horizontal dendriform structure associated to this $\quadri$ structure is
\begin{align*}
 T\west S  & =(T\ast\id)(\id\ast S)+T(\id\ast S\ast \id)\,,\\
 T\east S  & =(\id\ast T\ast
\id)S+(\id\ast T)(S\ast\id)\,. 
\end{align*}
This structure was found in~\cite[Corollary 4.14]{depaul} by other means.
The present work reveals that it is in fact the structure associated to a
(more fundamental) $\quadri$ structure.

The Baxter operators $\beta$ and $\gamma$ also give rise to dendriform structures on
$\End(A)$; these are obtained by the construction in~\ref{S:baxter-assoc}.
The dendriform structures are, respectively,
\[T\west S =T(\id\ast S)\,, \qquad T\east S=(\id\ast T)S\]
 and
\[ T\west S =T(S\ast\id)\,, \qquad T\east S=(T\ast\id)S\,. \]
The existence of these structures was announced in~\cite[Remark 4.15]{depaul}.

\section{The free $\quadri$ on one generator}\label{S:free}

\subsection{Free quadri-algebras}\label{S:free-quadri}
Let $V$ be a vector space. The {\it free quadri-algebra} $\calQ(V)$ on $V$ is a
$\quadri$ equipped with a map $i: V \to \calQ(V)$ which satisfies the following
universal property: for any linear map $f:V\to A$ where
$A$ is a $\quadri$, there is a unique $\quadri$ morphism $\Bar{f}:
\calQ(V)\to A$ such that
$\Bar{f} \circ i = f$. In other words, the functor $\calQ$ from vector spaces to
$\quadris$ is left adjoint to the forgetful functor.

The operad of $\quadris$ is a {\em non-$\Sigma$-operad}; in other words, the four operations of a
$\quadri$ have no symmetry, and the relations~\eqref{E:quadri} involve only monomials where
$x$, $y$ and $z$ stay in the same order. For this reason, the free $\quadri$
$\calQ(V)$ is of the form
$$\calQ(V) = \bigoplus_{n\geq 1} \calQ_n\ten V^{\ten n}\,.$$
Hence, $\calQ(V)$ is completely determined by the free  $\quadri$ on one generator 
$\calQ(\field)=\bigoplus_{n\geq 1}\calQ_n$.
Let $x$ be this generator. Then $x$ is a linear generator of $\calQ_1$ and the
elements  $x\se x$, $x \ne x$, $x \nw x$, $x \sw x$
form a basis  of $\calQ_2$. The space of four operations on three variables (with no
relations) is of dimension $2\times 4^2=32$. Since for $\quadris$ we have 9 linearly
independent relations, the space $\calQ_3$ is of dimension $32-9=23$.

\subsection{Conjecture}\label{S:conjecture} {\em The dimension of the vector space
$\calQ_n$  is}
\begin{equation}\label{E:dimension}
d_n := \frac{1}{n} \sum_{j=n}^{2n-1}{3n \choose n+1+j}{j-1 \choose j-n}\,.
\end{equation}
The first elements of this sequence are
\[ 1, 4, 23, 156, 1162, 9192, \ldots\]

According to~\cite[Theorem 2]{Fla}, $d_n$ is the number of {\em non-crossing connected
graphs} on $n+1$ vertices.

Let us give some evidence in favor of this conjecture.

\subsection{Dual quadri-algebras}\label{S:dual}
The operad $\calQ$ of $\quadris$ is binary (since it is generated by binary
operations) and is quadratic (since the relations involve only monomials with two operations).
Hence it has a dual operad $\calQ^!$(see~\cite{GK} for the seminal paper on the subject,
or~\cite[Appendix B]{L} for a short r\' esum\' e, or ~\cite{Fr} for a recent comprehensive
treatment).

Let us still denote the four dual operations by $\se$, $\ne$, $\nw$, $\sw$; these are linear
generators of
$\calQ^!_2$. The number of relations defining the binary quadratic operad $\calQ^!$
is 23 and so the dimension of $\calQ^!_3$ is $32-23=9$. Let us pick the following elements
as linear generators of $\calQ^!_3$:
\begin{align}
(x\nw y)\nw z & &(x\ne y)\nw z& &x\ne(y\ne z) \notag\\
(x\sw y)\nw z& &(x\se y)\nw z& &x\se(y\ne z)\label{E:gener-2}\\
(x\sw y)\sw z& &x\se (y\sw z)& &x\se(y\se z) \notag
\end{align}
From the 23 relations we deduce that any other monomial of degree 3 can be written as an algebraic
sum of these 9 elements.

\subsection{Proposition}\label{P:dimfree}
{\em The dimension of } $\calQ^!_n$ {\em is less than or equal to $n^2$.}
\bpf
Let $y$ be a monomial of degree $n$. It is determined by a planar binary tree (or parenthesizing)
where each node is labelled by one of the four generating operations. If, locally, a pattern
corresponding to one of the 23 discarded monomials of degree 3 appears, then we know that this
monomial can be rewritten in terms of other elements. Therefore, to generate $\calQ^!_n$
linearly, it  suffices to take the
$s_n$ monomials where only the 9 local patterns mentioned above appear.

Let $u_n$, $v_n$, $w_n$, $t_n$ be the number of these monomials whose lowest node is labelled respectively
by $\nw$, $\sw$, $\ne$ and $\se$. We have $s_n= u_n+ v_n+ w_n+ t_n$. From  the
choice of the 9 patterns it follows:
\begin{align*}
u_{n+1} &=  u_n+ v_n+ w_n+ t_n \\
v_{n+1} &=  v_n\\
w_{n+1} &= w_n\\
t_{n+1} &=  v_n+ w_n+ t_n 
\end{align*}
From these inductive relations (and $u_2=v_2= w_2=t_2=1$) it follows immediately that $s_n =
n^2$.
\epf

\subsection{Conjecture}\label{S:dimdual} ${\dim}\ \calQ^!_n=n^2$.

\subsection{Koszul duality}\label{S:Koszul} We further conjecture that the operad $\calQ$ is
 Koszul (cf. loc. cit.) Together with Conjecture~\ref{S:dimdual}, this implies
Conjecture~\ref{S:conjecture}. Indeed, they imply that the
generating series of
$\calQ$ and $\calQ^!$, that is, 
\[f(t) :=
\sum_{n\geq 1} (-1)^n \dim\calQ_n\, t^n \text{ \ and \ }g(t) :=
\sum_{n\geq 1} (-1)^n \dim\calQ_n^!\, t^n= \sum_{n\geq 1} (-1)^n n^2\, t^n = \frac{t(-1+t)}{(1+t)^3}\,,\]
 are inverse to each other with respect to composition:
$$ f(g(t))= t\,.$$
From here it follows that $\dim\calQ_n=d_n$ as in~\eqref{E:dimension};
see the Encyclopedia of Integer Sequences~\cite[A007297]{SlSeek}.

\subsection{Remark}\label{S:catalan} According to~\cite[equation (49)]{DB}, the
numbers $d_n$ satisfy the recursion
\[d_n=\sum_{i=1}^{n-1} d_id_{n-i}+\sum_{i,j=0}^{n-1}d_id_jd_{n-1-i-j}\,;\ \ d_0=1\,.\]
This is to be compared with the familiar recursion
 \[c_n=\sum_{i=0}^{n-1}c_ic_{n-1-i}\,;\ \ c_0=1\]
for the Catalan numbers $c_n$. It is known that $c_n$ is the dimension of the
homogeneous component of degree $n$ of the dendriform operad~\cite{L} and that the
inverse series  is $\sum_{n\geq 1} (-1)^n n\, t^n$ (compare with~\ref{S:Koszul}).

 \subsection{A quadri-algebra of shuffles}\label{S:shuffles}
 Recall the $\quadri$ $\field S_{\geq
2}=\bigoplus_{n\geq 2}\field S_n$ of Section~\ref{S:permutations}. Consider the subspace
 $Q_n$ of $\field S_{2n}$ spanned by the set
$\Sh(2,2,\ldots,2)$ of \odash{(2,2,\ldots,2)}{shuffles}. These are permutations $\zeta\in
S_{2n}$ of the form
\[\zeta=a_1b_1a_2b_2\ldots a_nb_n\]
with $a_i<b_i$ for every $i$. Let $Q:=\bigoplus_{n\geq 1} Q_n$. Clearly,
$Q$ is a $\quadrisub$ of $\field S_{\geq 2}$.
Moreover, $Q$ is graded if we declare that the elements of $Q_n$ have degree $n$.
Note that
\[\dim Q_n=\frac{(2n)!}{2^n}\,.\]

Let $F$ be the free $\quadri$ on one generator. Let $\iota:F\to Q$ be the unique
morphism of $\quadris$ that sends the generator to $12$ (the identity permutation in
$S_2$). 

\subsection{Conjecture}\label{C:shuffles}
{\em The map } $\iota:F\to Q$ {\em is injective. In other words, the $\quadrisub$ of $Q$
generated by $12$ is free.}

We discuss some evidence for this conjecture, which relates $F$ to the free dendriform
algebra and $Q$ to the dendriform algebra of permutations.

Let $Y$ be the free dendriform algebra on one generator~\cite{L}. Consider the
unique morphisms of dendriform algebras 
\[\psi_h:Y\to F_h \text{ \ and \ }\psi_v:Y\to F_v\]
 which send the generator of $Y$ to the generator of the $\quadri$ $F$.

Let $A:=\bigoplus_{n\geq 1}\field S_n$. There are two dendriform
structures on this space. With the notations of~\ref{S:permutations}, these structures are
\[\sigma \west \tau =
\sum_{\zeta\in\Sh(p,q),\zeta^{-1}(1)=1}\zeta\cdot(\sigma\times\tau)\,,
\qquad
\sigma \east \tau =
\sum_{\zeta\in\Sh(p,q),\zeta^{-1}(1)=p+1}\zeta\cdot(\sigma\times\tau)\]and
\[\sigma \north \tau =
\sum_{\zeta\in\Sh(p,q),\zeta^{-1}(p+q)=p}\zeta\cdot(\sigma\times\tau)\,,
\qquad
\sigma \south \tau =
\sum_{\zeta\in\Sh(p,q),\zeta^{-1}(p+q)=p+q}\zeta\cdot(\sigma\times\tau)\,.\]
We denote the first structure by $A_h$ and the second by $A_v$. They are
the horizontal and vertical dendriform algebra structures corresponding to the
$\quadri$ $\field S_{\geq 2}$ of~\ref{S:permutations}, enlarged by the component $\field S_1$ of degree
one.

Let
\[\alpha_h:Y\to A_h \text{ \ and \ } \alpha_v:Y\to A_v\]
be the unique morphisms of dendriform algebras 
 which send the generator of $Y$ to the permutation $1\in S_1$.

Consider now the $\quadri$ $Q$ of~\ref{S:shuffles} and the corresponding
dendriform algebras $Q_h$ and $Q_v$.

\subsection{Proposition}\label{P:twosquares}  {\em There are morphisms of dendriform algebras}
\[\varphi_h:A_h\to Q_h \text{ \ and \ } \varphi_v:A_v\to Q_v\]
{\em such that the following diagrams commute:}
\[\xymatrix{
F_h\ar[r]^{\iota} & Q_h\\
Y\ar[u]^{\psi_h}\ar[r]_{\alpha_h} & A_h\ar[u]_{\varphi_h} }
\qquad
\xymatrix{
F_v\ar[r]^{\iota} & Q_v\\
Y\ar[u]^{\psi_v}\ar[r]_{\alpha_v} & A_v\ar[u]_{\varphi_v} }
\]
\bpf Consider the maps $\Hat{\varphi}_h$ and $\Hat{\varphi}_v:\Sh(2,2,\ldots,2)\to S_n$
defined by
\[\Hat{\varphi}_h(a_1b_1a_2b_2\ldots a_nb_n)=\st(a_1a_2\ldots a_n) \text{ \ and \ }
\Hat{\varphi}_v(a_1b_1a_2b_2\ldots a_nb_n)=\st(b_1b_2\ldots b_n)\,,\]
where for any sequence of $n$ distinct integers $a_i$, $\st(a_1a_2\ldots a_n)$ denotes the
unique permutation $\sigma\in S_n$ such that
\[\sigma(i)<\sigma(j) \iff a_i<a_j\,.\]

Define $\varphi_h:A_h\to Q_h$ and $\varphi_v:A_v\to Q_v$ by
\[\varphi_h(\sigma)=\sum_{\Hat{\varphi}_h(\zeta)=\sigma}\zeta \text{ \ and \ }
\varphi_v(\sigma)=\sum_{\Hat{\varphi}_v(\zeta)=\sigma}\zeta\,. \]
It is easy to check that these are morphisms of dendriform algebras. 

Now, the composite $\varphi_h\alpha_h$ is a morphism of dendriform algebras which sends
the generator of $Y$ to $\varphi_h(1)=12$. Since the same is true of the composite
$\iota\psi_h$, the first square above must commute. The other commutativity is similar.
\epf

It is known that the maps $\alpha_h$ and $\alpha_v$ are injective. 
We view this as evidence in favor of Conjecture~\ref{C:shuffles}.

The maps $\alpha_h$ and $\alpha_v$ admit explicit combinatorial descriptions in terms of
planar binary trees. We expect similar descriptions for the maps $\psi_h$, $\psi_v$
and $\iota$.

\end{document}